\documentclass[11pt,reqno]{amsart}
\usepackage[T1]{fontenc}
\usepackage{lmodern,mathtools} 
\usepackage[margin=1.2in]{geometry}
\usepackage[dvipsnames]{xcolor}
\usepackage[colorlinks=true, pdfstartview=FitV, linkcolor= blue!80!black, citecolor=red!80!black, urlcolor=blue, pagebackref=false]{hyperref}
\newtheorem{theorem}{Theorem}[section]
\newtheorem{lemma}[theorem]{Lemma}
\newtheorem{corollary}[theorem]{Corollary}
\newtheorem{conjecture}[theorem]{Conjecture}
\theoremstyle{remark}
\newtheorem{remark}[theorem]{Remark}
\newcommand{\R}{\mathbb{R}}
\newcommand{\C}{\mathcal C}
\newcommand{\Haus}{\mathcal H}
\newcommand{\HS}{\mathrm{HS}}
\newcommand{\restriction}{\vert}
\newcommand{\e}{\varepsilon}
\numberwithin{equation}{section}

\DeclareMathOperator{\tr}{tr}
\DeclareMathOperator{\dist}{dist}

\title[Critical sets of harmonic functions]
{A near-quadratic bound for critical sets \\of harmonic functions}

\author[Z. Shen]{Zhongwei Shen}
\address{Zhongwei Shen: Institute for Theoretical Sciences, Westlake University, No. 600 Dunyu Road, Xihu District, Hangzhou,
Zhejiang 310030, P.R. China.}
\email{shenzhongwei@westlake.edu.cn}

\author[J. Zhuge]{Jinping Zhuge}
\address{Jinping Zhuge: Morningside Center of Mathematics, Academy of Mathematics and systems science, Chinese Academy of Sciences, Beijing 100190, China}
\email{jpzhuge@amss.ac.cn}

\thanks{J. Zhuge is partially supported by NNSF of China (Nos. 12494541, 12288201, 12471115).}

\begin{document}

\begin{abstract}
Let $u$ be a nonconstant harmonic function in $B_1 \subset \R^n$ with $n \ge 3$. Denote by $N$ the doubling index of $|\nabla u|$ in $B_1$. Let $\C := \{ x\in B_1: \nabla u(x) = 0 \}$ be the critical set. We show that 
$$\Haus^{n-2}(\C \cap B_{1/2}) \le C_n N^2 \log^2(N+2).$$ 
This estimate is sharp up to a squared logarithm loss.

\end{abstract}

\maketitle

\section{Introduction}

Let $u$ be a nonconstant harmonic function in $B_1 = B_1(0) \subset \R^n$ with $n \ge 2$.
This paper is concerned with quantitative estimate for the critical
set
\begin{equation*}
 \C=\{x\in B_1:\nabla u(x)=0\}.
\end{equation*}
Set
\begin{equation}\label{eq:normalization}
    N:= \log \frac{\sup_{B_1}| \nabla u|}{ \sup_{B_{1/2}}|\nabla u|}.
\end{equation}
Let $m(x)$ denote the vanishing order of $u - u(x)$ at $x\in B_1$.
Thus $m(x)-1$ is the vanishing order of $ \nabla u $ at $x$. For $x\in \C$, we have $2\le m(x) < \infty$. Without loss of generality, we assume $\C \cap B_{1/2}$ is not empty. This particularly rules out the trivial case that $u$ is an affine function. Under the assumption, $N \ge c_n$ for some dimensional constant $c_n > 0$.

The following is the main theorem of this paper.

\begin{theorem}\label{thm:main}
Let $n\geq3$ and let $u$ satisfy the assumptions above. Then
\begin{equation}\label{eq:main}
 \int_{\C\cap B_{1/2}}(m-1)\,d\Haus^{n-2}
 \leq C_n N^2\log^2(N+2),
\end{equation}
where $C_n$ depends only on $n$.
In particular,
\begin{equation}\label{eq:unweighted}
 \Haus^{n-2}(\C\cap B_{1/2})\leq C_n N^2\log^2(N+2).
\end{equation}
\end{theorem}

The proof of \eqref{eq:main} is based on the following identity, which is local and
holds on an arbitrary domain. 

\begin{theorem}\label{thm:identity}
Let $u$ be a nonconstant real harmonic function on a domain
$\Omega\subset\mathbb R^n$, $n\geq3$. Set $g = |\nabla u|$ and define
$F$ by the  Laplacian $\Delta \log g$ on $\{g>0\}$ and by zero
on $\C $. Then $\log g$ and $  F$ are locally integrable functions in $\Omega$, and
\begin{equation}\label{eq:distribution}
 \Delta \log g=F\,dx+2\pi(m-1)\,\Haus^{n-2}\restriction_{\C}
 \qquad\text{in }\mathcal D'(\Omega).
\end{equation}
On $\{g>0\}$, the density $F$ is given by
\begin{equation}\label{eq:regular-density}
 F=\frac{|\nabla^2u|_{\HS}^2-2|\nabla g|^2}{g^2},
\end{equation}
where $|\nabla^2 u|_{\HS}$ is the Hilbert-Schmidt norm of the Hessian $\nabla^2 u$ (see \eqref{eq.HS}).
\end{theorem}

For $n = 2$, the critical points are isolated and the number of critical points with multiplicity is bounded linearly by the frequency function (or doubling index). In fact, it has been shown in \cite{Han07,HLbook} that
\begin{equation*}
    \sum_{p\in \C \cap B_{1/2}} (m(p) - 1) \le CN.
\end{equation*}
A two-dimensional identity of \eqref{eq:distribution} also implies the above sharp result for harmonic functions as $F = 0$ in this case. An alternative proof was given in \cite{CZZ25} for planar elliptic equation with H\"{o}lder continuous coefficients.  For $n \ge 3$, Lin \cite[Theorem 2.1]{Lin} showed that the critical set of a harmonic function has Hausdorff dimension at most $n-2$ and conjectured
 \cite[Section 3, Conjecture 2]{Lin} an upper bound for the singular part of the nodal set, i.e., 
\begin{equation*}
    \Haus^{n-2}(\{ x\in B_{1/2}: |\nabla u(x)| = u(x) = 0 \}) \le C_n \Lambda^2,
\end{equation*}
where $\Lambda$ denotes the local Almgren frequency function, which is equivalent to the doubling index $N$ if $u(0) = 0$.
The conjecture for $n \ge 3$ was restated by Logunov and
Malinnikova \cite[Conjecture 1]{LM18} for the full critical set with the same quadratic upper bound in terms of the doubling index $N$ (known to be optimal for general harmonic polynomials), i.e.,
\begin{conjecture}[F. Lin]
For $n \ge 3$,
\begin{equation*}
    \Haus^{n-2}(\C \cap B_{1/2}) \le C_n N^2.
\end{equation*}
\end{conjecture}
Theorem~\ref{thm:main} gives this optimal upper bound up to a squared logarithm and includes the vanishing-order weight for any $n \ge 3$. The result may be extended to harmonic functions on real analytic manifolds.

Let us review the literature.
Hardt and Simon \cite[Eq.~(3.2) and Theorem~(3.3)]{HS89} obtained the optimal quadratic degree bounds for the critical sets of harmonic polynomials by algebraic slicing. Hardt, M. Hoffmann-Ostenhof, T. Hoffmann-Ostenhof and Nadirashvili \cite[Theorem~1.1 and Lemma~5.1]{HHHN} proved local finiteness and rectifiability of the critical sets for smooth elliptic equations. See also \cite{CNV15} for tubular volume estimates and quantitative stratification. Naber and Valtorta \cite[Theorem~1.1 and Remark~1.2]{NV17} then obtained an explicit upper bound for the $(n-2)$-dimensional Hausdorff measure in a form of $C^{\Lambda^2}$. Foster recently improved the bound to a quasipolynomial bound $\Lambda^{C\log \Lambda}$ in \cite{Fos25} for harmonic functions.

Now we briefly discuss the proof of the main theorem.
The pointwise estimates in Section~\ref{sec:linear-algebra} bound $q=(-F)_+$ by the two-dimensional Jacobian of $\nabla(\nabla u/|\nabla u|)$ and $q^{3/2}$ by the three-dimensional Jacobian of $\nabla^2 u$ divided by $g^3$. Section~\ref{sec:sublevel} gives exponential decay of gradient level sets and integrability of $\log |\nabla u|$. In Section~\ref{sec:poly}, we estimate the two Jacobians on slices of polynomial approximations by the area formula and B\'ezout's inequality. These estimates give the required integrability of $F$ and the bound for $q$. Theorems~\ref{thm:main} and~\ref{thm:identity} are proved in Section~\ref{sec:measure}.

The logarithmic gradient has a classical interpretation in the plane. If $u = u(x_1,x_2)$ is harmonic in two variables, the function $v = \partial_1 u-i \partial_2 u$ is holomorphic, and $\Delta \log |v|$ counts its zeros with multiplicity. Alessandrini \cite[Lemmas~2.1 and~2.3]{Alessandrini} used a corresponding logarithmic quantity for planar elliptic equations. See also Alessandrini and Magnanini \cite[Theorem~1.2, Eq.~(1.6)]{AM92} for an index formula, and Dong \cite[Theorem~3.3]{Dong} for a related distributional inequality arising from $\log(|\nabla u|^2 + \lambda u^2/2)$ at singular nodal points of surface eigenfunctions. An explicit distributional identity on a two-dimensional surface can also be found in \cite{CM01}.
In dimensions $n\geq3$, the analogue is established in Theorem \ref{thm:identity}, where the regular density $F$ in \eqref{eq:distribution} may take both signs. The main challenge is to bound its negative part $q = (-F)_+$ (the positive part can also be bounded; see Remark \ref{rmk.F}) and then identity the singular part $\Haus^{n-2}\restriction_{\C}$ by using basic tools from geometric measure theory.

\textbf{AI Use Statement.}
The proof of the main theorem was developed with the assistance of OpenAI’s ChatGPT. We believe that the approach was inspired in part by the work of Hardt and Simon [12], particularly in its use of Jacobians, gradient level sets, algebraic slicing arguments, the area formula, and zero counting for polynomials. Both authors have carefully verified the proof and substantially revised the argument, adding further details and references. They take full responsibility for the accuracy of the paper’s content.

\section{Pointwise Jacobian bounds}\label{sec:linear-algebra}

Let $u$ be a nonconstant harmonic function in $B_1=B(0, 1)$. We normalize $\| \nabla u \|_{L^\infty(B_1)} = 1$.
This leaves $N$, $\C$ and $m$ invariant. Let $G=\nabla u$ and $g=|\nabla u|$. Let $A=\nabla G=\nabla^2u$, and
$\nu=G/g$ where $g>0$. 
For an $n\times n$ matrix $B$, let $\{ \sigma_k(B) \}$ denote the singular values of $B$ in nonincreasing order. If $B$ is symmetric, then $\sigma_k(B) = |\lambda_k(B)|$, where $\{ \lambda_k(B) \}$ are eigenvalues of $B$. The $j$-dimensional Jacobian is defined by
\begin{equation}
    J_j(B) = \prod_{k=1}^j \sigma_k (B).
\end{equation}
Namely, $J_j(B)$ is the product of 
$j$ largest singular values of $B$.  We use $|\cdot|_{\HS}$
for the Hilbert-Schmidt norm. If $B = (b_{ij})$, then
\begin{equation}\label{eq.HS}
    |B|_{\HS}^2 = \sum_{i,j} |b_{ij}|^2 .
\end{equation}
Note that $|B|_{\HS}^2 = \tr(B^{\mathsf T} B) = \sum_{k=1}^n \sigma_k^2$.

\begin{lemma}\label{lem:pointwise}
In $\{g>0\}$, set $q=(-\Delta\log g)_+$. Then
\begin{equation}\label{eq:pointwise}
  q\leq C_nJ_2(\nabla \nu)\qquad
  \text{ and } \qquad
  q^{3/2}\leq C_n\frac{J_3(A)}{g^3}.
\end{equation}
\end{lemma}

\begin{proof}
Since the components of $G$ are harmonic, a direct computation shows that 
$$
\Delta(g^2)=2|A|_{\HS}^2 \quad \text{ and } \quad \nabla g=A\nu \quad \text{  on } \{g>0\}.
$$
Consequently,
\begin{equation}\label{eq:bochner}
 g^2\Delta\log g=|A|_{\HS}^2-2|A\nu|^2.
\end{equation}
Let $\lambda_i$ be the eigenvalues of $A$, ordered so that the singular value of $A$,
$\sigma_i=|\lambda_i|$, is decreasing. Since $\sum_i\lambda_i=0$ (due to $\Delta u = \tr(A) = 0$),
\begin{equation}
 \sigma_1\leq(n-1)\sigma_2\qquad
 \text{ and } \qquad
 \sigma_1-\sigma_2\leq|\lambda_1+\lambda_2|\leq(n-2)\sigma_3.
\end{equation}
Using $|A\nu|\leq \sigma_1$ (since $|\nu| = 1$) and $|A|_{\HS} \ge \sigma_1^2 + \sigma_2^2$ in \eqref{eq:bochner} yields
\begin{equation}\label{eq.g2q}
  g^2q\leq (  2|A\nu|^2 - |A|_{\HS}^2)_+ \leq \sigma_1^2-\sigma_2^2\leq n(n-2)\sigma_2\sigma_3.
\end{equation}
By the identity $g \nabla \nu=(I-\nu\nu^{\mathsf T})A$, we have
\begin{equation*}
 g^2 (\nabla \nu)^{\mathsf T}(\nabla \nu)=A^2-(A\nu)(A\nu)^{\mathsf T}.
\end{equation*}
For $k=1,2$, the span of the first $k+1$ eigenvectors of $A^2$
intersects $(A\nu)^\perp$ in a subspace of dimension at least $k$. For any vector $\xi$ on that subspace, we have
\begin{equation*}
    \xi^{\mathsf T} \cdot A^2 \xi = \xi^{\mathsf T} (A^2-(A\nu)(A\nu)^{\mathsf T}) \xi.
\end{equation*}
The min-max principle therefore bounds the two largest eigenvalues of $A^2-(A\nu)(A\nu)^{\mathsf T}$ from below by the second and third eigenvalues of $A^2$, $\sigma_2^2$ and $\sigma_3^2$, respectively.
It follows that $g^2 J_2(\nabla \nu)\geq \sigma_2\sigma_3$. This, together with \eqref{eq.g2q}, implies the first inequality in \eqref{eq:pointwise}. Obviously, $(\sigma_2 \sigma_3)^{3/2} \le \sigma_1 \sigma_2 \sigma_3$. Thus, \eqref{eq.g2q} also yields the second inequality in \eqref{eq:pointwise}.
\end{proof}


In the next lemma, we show that $J_j(B)$ is bounded (actually comparable) by the sum of the absolute values of its $j$-by-$j$ minors. Let $I = \{ i_1 < i_2 < \cdots < i_j \}$ and $K = \{ k_1 < k_2 < \cdots < k_j \}$ be subsets of $\{ 1,2,\ldots, n\}$. Let $\mathcal{S}_j$ be the set of all such index subsets of $\{1,2, \ldots, n \}$ with $j$ elements. Let $B_{I, K}$ be the $j$-by-$j$ minor taken from $B$ with $B_{I,K} = (b_{ik})_{i\in I, k\in K}$.

\begin{lemma}\label{lem:Jacobian-det}
Let $B$ be an $n\times n$ matrix.
    For $1\le j\le n$, it holds
    \begin{equation}
        J_j(B) \le \sum_{I, K \in \mathcal{S}_j} |\det (B_{I,K})|.
    \end{equation}
\end{lemma}

\begin{proof}
    Let $\sigma_1\ge\cdots\ge \sigma_n $ be the singular values of $B$, and let $U, V \in \R^{n\times j}$ have as columns the corresponding first $j$ left and right singular vectors, respectively. Thus, the columns of $U$ and $V$ are orthonormal vectors in $\R^n$. Then, by the singular value decomposition,
    \begin{equation*}
    U^{\mathsf T}BV=\operatorname{diag}(\sigma_1,\dots,\sigma_j)
    \qquad \text{ and } \qquad
    \det(U^{\mathsf T}BV) = J_j(B).
    \end{equation*}
    Applying Cauchy-Binet formula twice gives
    \begin{equation}\label{eq.JjB=}
    J_j(B)=\sum_{I,K\in\mathcal{S}_j}
    \det(U_I)\,\det(B_{I,K})\,\det(V_K),
    \end{equation}
where \(U_I\) and \(V_K\) are the submatrices formed by taking rows \(I\) and \(K\), respectively. Each column of \(U_I\) or \(V_K\) has norm at most one. Thus Hadamard’s inequality implies \(|\det U_I|,|\det V_K|\le1\). Finally, taking absolute values of  \eqref{eq.JjB=} yields the desired estimate.
\end{proof}

\section{Level set estimates}\label{sec:sublevel}

Recall that $g=|\nabla u|$, where $u$ is a nonconstant harmonic function in $B_1$. A volume estimate of the form below follows by applying the level set estimate of Logunov and Malinnikova
\cite[Lemma~5.3]{LM18} to the gradient in interior cubes (together with a covering argument). We provide a different proof by using holomorphic extension \cite{Lin}.

\begin{lemma}\label{lem:sublevel}
For $t\geq 0$,
\begin{equation}\label{eq:sublevel}
 |\{x\in B_{3/4}:g(x)<e^{-t}\}|\leq C_n e^{-c_nt/N},
\end{equation}
and
\begin{equation}\label{eq.log-g-L1}
    \int_{B_{3/4}}|\log g|\,dx\leq C_n N .
\end{equation}
\end{lemma}

\begin{proof}
    We will prove that  for some $c_n > 0$ and $p = c_n/N$, 
    \begin{equation}\label{eq.g^-p}
        \int_{B_{3/4}} g^{-p} dx \le C_n.
    \end{equation}
    We first reduce the gradient estimate to an estimate for a scale harmonic function.
    By the doubling condition, there exists $x_* \in B_{1/2}$ such that $g(x_*) \ge \frac12 e^{-N}$. Fix such $x_*$ and set
    \begin{equation*}
        f(x) =\frac{G(x_*)}{g(x_*)} \cdot G(x).
    \end{equation*}
    Note that $f$ is harmonic in $B_1$ and $|f| \le g \le 1$ and $f(x_*) = g(x_*) \ge \frac12 e^{-N}$. Thus to prove \eqref{eq.g^-p}, we only need to prove the corresponding estimate for $f$.

    We will use the method of holomorphic extension. For $\theta \in \mathbb{S}^{n-1}$, let $x_* + R_\theta \theta \in \partial B_{3/4}$. Define $f_\theta(s) = f(x_* + s R_\theta \theta)$. The Poisson formula on balls of radius $1/8$ centered in $\overline{B_{3/4}}$ implies that $f_\theta(s)$ can be extended holomorphically to a complex function $f_\theta(z)$. Moreover, there exist $\delta_n > 0$ and $M_n > 0$, independent of $\theta$, such that
    \begin{equation}
        |f_\theta(z)| \le M_n, \quad \text{for all } z\in \Omega:= \{ z\in \mathbb{C}: \dist(z, [0,1]) < \delta \}.
    \end{equation}
    Fix a conformal mapping $\phi: \Omega \mapsto \mathbb{D}$, where $\mathbb{D}$ is the unit disk in $\mathbb{C}$, with $\phi(0) = 0$ and $\phi'(0) > 0$. By symmetry of $\Omega$, $\phi([0,1]) = [0, a]$, $a := \phi(1) < 1$ and $\sup_{[0,a]} |(\phi^{-1})' | \le C_n$.
    Put
    \begin{equation}
    F(z) = \frac{f_\theta\circ\phi^{-1}(z)} {M_n}, \qquad z\in \mathbb{D}.
    \end{equation}
    Then $\|F\|_{L^\infty(\mathbb D)}\le1$ and
    \begin{equation*}
        N_1 := \log\frac1{|F(0)|}\le N+\log(2M_n)\le C_n N. 
    \end{equation*}
    Fix $b = (1+a)/2 < 1$.
    By \cite[Lemma 2.3.2]{HLbook}, we know the number of zeros of $F$ in $\{ |z| < b \}$ is bounded by $C_n N$. Suppose $\alpha_1, \alpha_2,\dots, \alpha_k$ are the zeros of $F$ in $\{ |z| \le b \}$, with multiplicity. Let $B_\alpha = (z-\alpha)/(1-\bar{\alpha} z)$ and
    \begin{equation}
        Q(z) = F(z) / \prod_{j=1}^k B_{\alpha_j}(z).
    \end{equation}
    Thus $Q(z)$ is holomorphic in $\mathbb{D}$ with no zeros in $\{ |z| < b \}$. Since $|B_\alpha| = 1$ on $\partial \mathbb{D}$, we have $|Q(z)| = |F(z)| \le 1$ on $\partial \mathbb{D}$. The maximum modulus principle implies $|Q| \le 1$ in $\mathbb{D}$. Thus $-\log |Q|$ is nonnegative and harmonic in ${|z| < b}$. Since $|Q(0)| \ge |F(0)| = |f(x_*)| \ge \frac12 e^{-N}$, the Harnack's inequality implies
    \begin{equation*}
        -\log |Q(z)| \le -C\log |Q(0)| \le C_n N_1 \le C_n N, \quad \text{for all } |z| \le a.
    \end{equation*}

    Now, choose $c_n > 0$ such that $p = c_n /N$ and $kp \le 1/2$, where $k$ is the number of zeros of $F$ in $\{|z| < b \}$. The arithmetic-geometric mean inequality and $|B_\alpha| \le 1$ yield
    \begin{align*}
    |F(t)|^{-p}
    &\le e^{C_n Lp}\,\frac1k\sum_{j=1}^k |B_{\alpha_j}(t)|^{-kp}\\
    &\le \frac{C_n}{k}\sum_{j=1}^k |t-\operatorname{Re}\alpha_j|^{-1/2}.
\end{align*}
Since $\int_0^a |t-\alpha|^{-1/2} dt \le C$ for any $\alpha \in \R$, we obtain
\begin{equation}
    \int_0^a |F(t)|^{-p} \, dt \le C_n.
\end{equation}
By a change of variables, we have
    \begin{equation*}
    \int_0^1 |f_\theta(s)|^{-p}\,ds
    =M_n^{-p}\int_0^a |F(t)|^{-p}|(\phi^{-1})'(t)|\,dt\le C_n.
    \end{equation*}
    This estimate is uniform for $\theta$. Since $g \ge |f|$, by spherical coordinates,
    \begin{equation*}
    \int_{B_{3/4}} g^{-p}\,dx
    \le \int_{\mathbb S^{n-1}}R_\theta^n
         \int_0^1 |f_\theta(s)|^{-p}s^{n-1}\,ds\,d\theta\le C_n .
    \end{equation*}
    Therefore, for $t\ge0$,
    \begin{equation*}
     |\{x\in B_{3/4}: g(x)<e^{-t}\}|
    \le e^{-pt}\int_{B_{3/4}} g^{-p}\,dx\le C_n e^{-c_nt/N}.
    \end{equation*}
    Finally, since $0\le g\le1$,  the layer-cake formula gives
    $$
    \aligned
    \int_{B_{3/4}} |\log g|\,dx
     & =\int_0^\infty |\{x\in B_{3/4}:g(x)<e^{-t}\}|\,dt\\
    & \le C_n\int_0^\infty e^{-c_nt/N}\,dt\le C_n N.
    \endaligned
    $$
    The proof is complete.
\end{proof}



\section{Polynomial approximation and Jacobian integrability}\label{sec:poly}

Let $G=\nabla u$. For $T\geq1$ and fixed $b>0$, there is a real vector polynomial $P$
of degree $d\leq C_{n,b}T$ such that
\begin{equation}\label{eq:approx}
 \|P-G\|_{C^1(B_{3/4})}\leq e^{-bT}.
\end{equation}
This property follows from the harmonicity of $G$ and is independent of the doubling index.
Indeed, we can expand the components of $G$ in homogeneous harmonic polynomials about the origin,
\begin{equation*}
    G = \sum_{\ell \ge 0} G_\ell,
\end{equation*}
where $G_\ell$ are vector-valued homogeneous harmonic polynomials of degree $\ell$. By orthogonality of spherical harmonics, and our normalization $|G| \le 1$,
\begin{equation*}
    \int_{\partial B_1} |G|^2 = \sum_{\ell \ge 0} \int_{\partial B_1} |G_\ell|^2 \le C_n.
\end{equation*}
Note that for each $x \in B_1$, we have $G_\ell(x) = |x|^\ell G_\ell(x/|x|)$. Now we set
\begin{equation}
    P = \sum_{\ell \ge 0}^d G_\ell.
\end{equation}
Hence, by homogeneity and orthogonality,
\begin{equation}
    \int_{\partial B_{7/8}} |G - P|^2 = C\sum_{\ell > d} \Big( \frac{7}{8} \Big)^{2\ell} \int_{\partial B_1} |G_\ell|^2 \le C_n \Big( \frac{7}{8} \Big)^{2d}.
\end{equation}
Finally, since $G-P$ is harmonic, interior estimates imply \eqref{eq:approx} by choosing $C_{n,b}$ large enough depending on $n$ and $b$.
This also gives $|\nabla G|+|\nabla P|\leq C_n$ on $B_{3/4}$.

We will need the area formula and B\'{e}zout inequality stated in the following two lemmas.
\begin{lemma}[Area formula]\label{lem:area}
Let $F:U\to\mathbb R^k$ be a $C^1$ map on an open set
$U\subset\mathbb R^k$, and let $E\subset U$ be Lebesgue measurable.
Then
\begin{equation}\label{eq:area}
 \int_E|\det \nabla F(x)|\,dx
 =\int_{\mathbb R^k}
 \#\{x\in E: F(x)=y \text{ and } \det \nabla F(x)\ne0\}\,dy .
\end{equation}
The equality also holds with all pre-images counted on the
right-hand side.
\end{lemma}

This is the equal-dimensional case of the area formula
\cite[Chapter~2, Theorem~3.3]{SimonGMT}. We add the constraint $\det \nabla F(x) \neq 0$ on the right-hand side because we can apply the formula to $\widetilde{E} = E \setminus \{\det \nabla F = 0 \}$, which is also Lebesgue measurable. The integral on the left-hand side of \eqref{eq:area} does not distinguish between $E$ and $\widetilde{E}$.

\begin{lemma}[B\'ezout's inequality]\label{lem:bezout}
Let $H_1,\ldots,H_k$ be nonzero complex polynomials in $k$
variables of positive degrees $d_1,\ldots,d_k$. The number
of isolated common zeros in $\mathbb C^k$ is at most
$\prod_{i=1}^k d_i$. In particular, if $H=(H_1,\ldots,H_k)$ is
a real polynomial map of degree at most $d$, then for every
$y\in\mathbb R^k$,
\begin{equation}\label{eq:bezout}
 \#\{x\in\mathbb R^k: \ H(x)=y \text{ and } \det \nabla H(x)\ne0\}\leq d^k .
\end{equation}
\end{lemma}

The proof of this lemma can be found in e.g. \cite[Theorem~1.1]{Sot}.

Recall that $g=|\nabla u|$ and $q = (-\Delta \log g)_+$. The next two lemmas are the key ingredients for the estimates of $q$, which are independent of the doubling index.

\begin{lemma}\label{lem:two}
For every $T\geq1$,
\begin{equation}\label{eq:two}
 \int_{B_{3/4} \cap\{g\geq e^{-T}\}}q\,dx\leq C_nT^2.
\end{equation}
\end{lemma}

\begin{proof}
Set $\delta=e^{-T}$. Choose $P$ by \eqref{eq:approx} so that 
$$
\|P-G\|_{C^1(B_{3/4})}\leq \delta^3 = e^{-3T}
$$
and $d\leq C_nT$.
On $\{g\geq\delta\}$, the segment from $G(x)$ to $P(x)$, 
$$
\big\{Z(t) = t G(x) + (1-t) P(x): t\in [0,1] \big\}
$$
stays $\delta/2$-away from the origin, i.e., $|Z(t)| \ge \delta/2$ for $t\in [0,1]$. Let $\Phi(Z) := \widehat{Z} = Z/|Z|$. Then an algebraic computation shows that 
\begin{equation}\label{eq.DPhi}
    \nabla_Z \Phi(Z)=|Z|^{-1}(I-\widehat Z\otimes\widehat Z)
    \qquad \text{ and } \qquad |\nabla^2_Z\Phi(Z)|\leq C_n |Z|^{-2}.
\end{equation}
Write $\omega = \Phi(P) = P/|P|$ and recall $\nu = \Phi(G) = G/|G|$. The chain rule implies that 
$$
\nabla \omega = \nabla_Z \Phi(P) \nabla P\qquad
\text{ and } \qquad\nabla \nu = \nabla_Z \Phi(G) \nabla G.
$$
It follows that
    \begin{equation*}
    |\nabla \nu|+|\nabla \omega|\leq C_n\delta^{-1},
    \end{equation*}
and
 \begin{equation*}
    |\nabla \nu - \nabla \omega| \le C_n\delta^{-2}|G-P|+C_n\delta^{-1}|\nabla G- \nabla P|
 \le C_n\delta,
\end{equation*}
where we have used mean value theorem, \eqref{eq.DPhi}, and the fact $|\nabla P| + |\nabla G| \le C_n$.
By Lemma \ref{lem:Jacobian-det},
\begin{equation}
\begin{aligned}
    J_2(\nabla \nu) & \le \sum_{I,K\in \mathcal{S}_2} |\det (\nabla \nu)_{I,K}| \\
    & \le \sum_{I,K\in \mathcal{S}_2} |\det (\nabla \omega)_{I,K}| + C_n (|\nabla \nu| + |\nabla \omega|) |\nabla \nu-\nabla\omega| \\
    & \le \sum_{I,K\in \mathcal{S}_2} |\det (\nabla \omega)_{I,K}| + C_n.
\end{aligned}
\end{equation}
Consequently, it suffices to estimate $$
\int_{B_{3/4} \cap \{ g > \delta \}} |\det (\nabla \omega)_{I,K}|.
$$ Since the estimates are the same for every pair $I, K \in \mathcal{S}_2$, without loss of generality, we will only consider $I = K = \{1,2 \}$ below.

For $I = K = \{ 1,2\}$, we have $ $
\begin{equation}
    (\nabla \omega)_{I,K} = \begin{bmatrix}
        \partial_1 \omega_1 & \partial_2 \omega_1 \\
        \partial_1 \omega_2 & \partial_2 \omega_2
    \end{bmatrix}.
\end{equation}
Write $x = (x_1, x_2, x')$ and fix $x' \in \R^{n-2}$. We view
\begin{equation}\label{eq.defF}
    F(x_1, x_2) = (\omega_1, \omega_2) = \Big( \frac{P_1}{|P|}, \frac{P_2}{|P|} \Big)\qquad \text{ and }\qquad \nabla F(x_1,x_2) = (\nabla \omega)_{I,K}
\end{equation}
as functions of $(x_1,x_2)$.
Here and after, for simplicity, we omit the implicit dependence on $x'$.
Since $|P|$ does not vanish on $D_\delta:= B_{3/4} \cap \{ g > \delta \}$, $F$ is a smooth function in $D_\delta$ with range contained in $B_1 \cap \R^2$. Let 
$$
D_\delta(x') = \{(x_1,x_2): \ (x_1,x_2,x') \in D_\delta \}.
$$
It follows from the area formula in Lemma \ref{lem:area} with $k = 2$ that 
\begin{equation}\label{eq.DetDw}
\begin{aligned}
    & \int_{D_\delta(x')} |\det (\nabla \omega)_{I,K}|\, dx_1 dx_2 \\
    & = \int_{B_1\cap \R^2} \# \Big\{(x_1, x_2) \in D_\delta(x'): F(x_1, x_2) = (s,t)\text{ and } \det \nabla F(x_1,x_2) \neq 0 \Big\} \, ds dt.
\end{aligned}
\end{equation}
To count the number of solutions of $F(x_1,x_2) = (s,t)$ for each pair $(s,t)$, we use \eqref{eq.defF} to write it as
\begin{equation}
    H_1:= P_1^2 - s^2 |P|^2 = 0\qquad \text{ and } \qquad H_2:= P_2^2 - t^2 |P|^2 = 0.
\end{equation}
This is a system of polynomial equations in two variables with degrees at most $2d$. Write $H = (H_1,H_2)$. At any points satisfying $F(x_1,x_2) = (s,t)$, a straightforward calculation shows 
\begin{equation*}
    \nabla H(x_1,x_2) =
    \begin{bmatrix}
        2s|P|^2 & 0 \\
        0 & 2t|P|^2
    \end{bmatrix}
    \nabla F(x_1,x_2).
\end{equation*}
This shows $\det \nabla H(x_1,x_2) = 4st |P|^4 \det \nabla F(x_1,x_2)$ if $F(x_1,x_2) = (s,t)$. Hence, if $st \neq 0$, $\det \nabla H \neq 0$ if and only if $\det \nabla F(x_1,x_2) \neq 0$ at isolated zeros. Thus, by B\'{e}zout's inequality in Lemma \ref{lem:bezout}, if $st \neq 0$,
\begin{equation}
    \# \Big\{(x_1, x_2) \in D_\delta(x'): F(x_1, x_2) = (s,t) \text{ and } \det \nabla F(x_1,x_2) \neq 0 \Big\} \le 4d^2.
\end{equation}
Since the coordinate axes $st = 0$ have measure zero in $\R^2$, the right-hand side of \eqref{eq.DetDw} does not see this case. Consequently, we obtain
\begin{equation}
    \int_{D_\delta(x')} |\det (\nabla \omega)_{I,K}|\, dx_1 dx_2 \le Cd^2.
\end{equation}
This estimate is uniform in $x'$. Integrating over $x'$, we then obtain
\begin{equation}
    \int_{D_\delta} |\det (\nabla \omega)_{I,K}|\, dx \le Cd^2.
\end{equation}
Summing over all $I, K\in \mathcal{S}_2$, we get
\begin{equation}
    \int_{D_\delta} J_2(\nabla \nu) \le CT^2,
\end{equation}
which leads to \eqref{eq:two} by the first inequality in \eqref{eq:pointwise}.
\end{proof}

\begin{lemma}\label{lem:three}
For each integer $j\geq1$, let
$E_j = B_{3/4} \cap\big\{e^{-j-1}\leq g<e^{-j}\big\}$. Then
\begin{equation}\label{eq:three}
 \int_{E_j}q^{3/2}\,dx\leq C_nj^3.
\end{equation}
\end{lemma}

\begin{proof}
Set $\delta=e^{-j}$. In $E_j$, the second
inequality of \eqref{eq:pointwise} bounds $q^{3/2}$ by
$C_n\delta^{-3}J_3(\nabla G)$ in $E_j$. We need to extract a factor $\delta^3$ from $J_3(\nabla G)$. The proof is similar to Lemma \ref{lem:two}.

Let $P$ be a polynomial vector given by \eqref{eq:approx} such that $\| P - G \|_{C^1(B_{3/4})} \le \delta^3$ and the degree of $P$ satisfies $d \le Cj$.
Note that the derivatives $\nabla G, \nabla P$ are uniformly bounded on the patch in $B_{3/4}$. By Lemma \ref{lem:Jacobian-det},
\begin{equation}\label{eq.J3DG}
\begin{aligned}
    \int_{E_j} J_3(\nabla G) & \le \sum_{I, K\in \mathcal{S}_3 } \int_{E_j} |\det (\nabla G)_{I,K}| \\
    & \le \sum_{I, K\in \mathcal{S}_3 } \int_{E_j} |\det (\nabla P)_{I,K}| + C\delta^3.
\end{aligned}
\end{equation}

As before, it suffices to consider $I = K = \{ 1,2,3\}$. Write $z = (x_1,x_2,x_3)$ and $x = (z, x'')$, where $x''$ denotes the remaining variables. Put
\begin{equation}
    Q(z) = (P_1, P_2, P_3)(z),
\end{equation}
and for a fixed $x''$,
\begin{equation}
    E_j(x'') =\big \{ z: (z,x'') \in E_j \big\}.
\end{equation}
Thus, $\nabla_z Q = (\nabla P)_{I,K}$. Consequently, as in Lemma \ref{lem:area} by the area formula and B\'{e}zout inequality with $k = 3$, we arrive at
\begin{equation*}
\begin{aligned}
    \int_{E_j(x'')} |\det (\nabla P)_{I,K}| & = \int_{\{|z| \le 2\delta \}} \#\Big\{ z\in E_j(x''):\  Q(z) = y
    \text{ and }
    \det \nabla_z Q \neq 0 \Big\} \, dy\\  
    &\le C\delta^3 d^3.
\end{aligned}
\end{equation*}
Integrating over $x''$ and summing over $I,K\in \mathcal{S}_3$ in \eqref{eq.J3DG}, we get
\begin{equation}
    \int_{E_j} J_3(\nabla G) \le C\delta^3 j^3,
\end{equation}
which yields the desired estimate.
\end{proof}

\begin{corollary}\label{cor:q}
    The function $q$ is integrable in $B_{3/4}\setminus \C$ and
    \begin{equation}\label{eq:qmass}
        \int_{B_{3/4}\setminus \C }q\,dx\leq C_n N^2\log^2(N+2).
    \end{equation}
\end{corollary}
\begin{proof}
    By H\"older's inequality, \eqref{eq:sublevel}, and
\eqref{eq:three},
\begin{equation}\label{eq:shellmass}
  \int_{E_j}q\,dx
   \leq \left(\int_{E_j}q^{3/2}\,dx\right)^{2/3}|E_j|^{1/3}
   \leq C_nj^2 e^{-c_nj/N}.
\end{equation}
For any integer $J\geq 1$ (to be determined), split the domain $B_{3/4}\setminus \C$ into $\{ g \ge e^{-J} \} \cap B_{3/4}$ and $\cup_{j \ge J} E_j$. Then, by Lemma \ref{lem:two} and Lemma \ref{lem:three}, we have
\begin{equation}
\begin{aligned}
    \int_{B_{3/4}} q & = \int_{\{ g \ge e^{-J} \} \cap B_{3/4}} q + \sum_{j\geq J} \int_{E_j} q \\
    &\le C_nJ^2 + \sum_{j\geq J} j^2 e^{-c_nj/N} \\
    & \le C_nJ^2 + C_n N^3 e^{-c_n J/(2N)}.
\end{aligned}
\end{equation}
Choose $J=\lceil C_n N\log(N+2)\rceil$ with the coefficient large enough so that the second term is at most $C_nN ^2$. Then the first is bounded by $C_n N^2 \log^2 (N+2)$. This gives the desired estimate.
\end{proof}

\section{The distributional identity and the measure estimate}\label{sec:measure}

In this section, we prove Theorem \ref{thm:identity} and then show it implies Theorem \ref{thm:main}.

\begin{proof}[Proof of Theorem \ref{thm:identity}]
Let $B$ be any ball in $\Omega$ and we only need to prove the distributional identity \eqref{eq:distribution} in $(1/2)B$. Without loss of generality, assume $B = B_1$. Suppose $B_{1/2} \cap \C$ is not empty; for otherwise the identity is a direct computation since $\log g$ is smooth. Let $N$ be the doubling index given by \eqref{eq:normalization}.

Let $w = \log g = \log |\nabla u|$.
Set $w_\varepsilon=\tfrac12\log(g^2+\varepsilon^2)$, and let $F$
be the same as in Theorem~\ref{thm:identity}. 
The chain rule gives
\begin{equation}\label{eq:reg}
 \Delta w_\varepsilon=
 \left(\frac{g^2}{g^2+\varepsilon^2}\right)^2 F
 +\rho_\varepsilon\qquad
 \text{ and } \qquad
\rho_\varepsilon=\frac{\varepsilon^2|\nabla^2u|_{\HS}^2}
 {(g^2+\varepsilon^2)^2}
\end{equation}
almost everywhere. Indeed,
$\Delta(g^2)=2|D^2u|_{\HS}^2$ and
$\nabla(g^2)=2g\nabla g$ off $\C$, while $\C$ has Lebesgue
measure zero by \cite[Theorem 2.1]{Lin}. Since $q=(-F)_+=(F)_-$, the measures
$(\Delta w_\varepsilon+q)\,dx$ are nonnegative. For $0<\varepsilon\leq1$,
$w\leq w_\varepsilon\leq \tfrac12\log 2$ under the normalization
$\| \nabla u \|_{L^\infty(B_1)} = 1$. By Lemma~\ref{lem:sublevel} and monotonicity, $w_\varepsilon\to w$ in
$L^1(B_{3/4})$. It follows that
\begin{equation}\label{eq.mu}
 \mu=\Delta w+q\,dx
\end{equation}
is a positive distribution and hence a locally finite Radon measure~\cite[Theorem 6.22]{LL01}. In $B_1\setminus \C$, it holds $\mu = F_+\,dx$. Its local finiteness 
and \eqref{eq:qmass} give $F_+\in L^1_{\rm loc}(B_{3/4})$. Indeed, let $\eta \in C_0^\infty(B_{3/4})$ such that $\eta = 1$ in $B_{a}$ with $a\in (0,3/4)$ and $0\le \eta\le 1$. Then by positivity of $\mu$ and $\mu\restriction_{B_1\setminus \C} = F_+ dx$, we arrive at
\begin{equation}\label{eq.F+L1}
    \int_{B_{a}} F_+ dx \le \mu(\eta) = \int_{B_{3/4}} (w \Delta \eta + q \eta ) dx \le CN^2\log^2(N+2),
\end{equation}
where we have used Lemma \ref{lem:sublevel} and \eqref{eq:qmass} in the last inequality.
Consequently, we have
\begin{equation}\label{eq:measure-decomposition}
 \Delta w=F\,dx+\mu\restriction_{\C},
\end{equation}
where $F \in L^1_{\rm loc}(B_{3/4})$.
In addition, $\mu = F_+\, dx + \mu \restriction_\C$.

Fix $p\in\C \cap B_{1/2}$ and let $k=m(p)-1$. With
$M_p(r)=\sup_{B_r(p)}g$, the homogeneous Taylor expansion of $G=\nabla u$ gives $M_p(r)=c_p r^k(1+O_p(r))$ with $c_p>0$. 
As a result, the doubling index $\log(M_p(4r)/M_p(2r))$ tends to $k\log 2$ as $r \to 0$.
We will apply \eqref{eq:sublevel} and \eqref{eq:qmass} to the harmonic function
\begin{equation*}
 U_r(z) = \frac{u(p+4rz)-u(p)}{4rM_p(4r)}.
\end{equation*}
Observe that $\| \nabla U_r \|_{L^\infty(B_1)} = 1$ and
\begin{equation}\label{eq.Ur-log}
    \log |\nabla U_r(z)| = \log |(\nabla u)(p + 4rz)| - \log M_p(4r).
\end{equation}
Consequently,
\begin{equation}\label{eq.Ur-Deltalog}
    (-\Delta \log |\nabla U_r(z)|)_+ = (4r)^2 q(p + 4rz).
\end{equation}
It follows by applying \eqref{eq:sublevel}  to \eqref{eq.Ur-log} that
\begin{equation}\label{eq:local-growth-log}
\begin{aligned}
    \int_{B_{2r}(p)}|w-k\log r|\,dx &\le (4r)^n \int_{B_{1/2}}|\log |\nabla U_r(z)||\,dz + |B_{2r}||\log M_p(4r) - k \log r| \\
    &\le C_p r^n
\end{aligned}
\end{equation}
for all sufficiently small $r$, where $C_p$ is a constant depending on the doubling index at $p$ (independent of $r$). Moreover, by \eqref{eq.Ur-Deltalog} and \eqref{eq:qmass},
\begin{equation}
    \int_{B_{2r}(p)}q\,dx\le C_p r^{n-2}.
\end{equation}

Let $\chi_r$ be a nonnegative cutoff equal to one on $B_r(p)$,
supported in $B_{2r}(p)$, with $\chi_r\leq1$ and
$|\Delta\chi_r|\leq Cr^{-2}$. Since $\int\Delta\chi_r=0$,
\eqref{eq:local-growth-log} and positivity of $\mu$ imply
\begin{equation}\label{eq:local-growth-measure}
 \mu(B_r(p))\le
 \int(w-k\log r)\Delta\chi_r\,dx+\int\chi_rq\,dx
 \le C_p r^{n-2}.
\end{equation}
Thus $\mu\restriction_{\C}\ll\sigma$, where
$\sigma=\Haus^{n-2}\restriction_{\C}$. To see this, suppose that $E\subset \C\cap B_{1/2}$ has $\mathcal H^{n-2}(E)=0$. For each integer $j\ge1$, let $E_j$ consist of the points in $E$ where the estimate \eqref{eq:local-growth-measure} holds with $C_p=j$ and $r_p=1/j$. Then $E=\bigcup_{j\ge1}E_j$. Since $\Haus^{n-2}(E_j)=0$, for every $\e > 0$, we can cover $E_j$ by balls $B_{r_i}(p_i)$ centered in $E_j$ with $r_i<1/j$ and $\sum_i r_i^{n-2}<\varepsilon$. The uniform bound on $E_j$ gives
\begin{equation*}
\mu^*(E_j)
\le \sum_i\mu(B_{r_i}(p_i))
\le j\sum_i r_i^{n-2}
<j\e,
\end{equation*}
where $\mu^*$ is the outer measure. Since $\e$ is arbitrary, it follows that $\mu^*(E_j)=0$. Countable subadditivity then gives $\mu^*(E)=0$, which proves the absolute continuity.

By Radon-Nikodym theorem, we can write $\mu\restriction_{\C}=f\sigma$. Let
$\omega_{n-2}$ denote the volume of the unit ball in
$\mathbb R^{n-2}$. By \cite[Theorem~1.1 and Lemma~5.1]{HHHN},
$\sigma$ is locally finite and $\C$ is countably
$(n-2)$-rectifiable. By the density theorem~\cite[Theorem 3.2.19]{Federer} and the
differentiation of Radon measures~\cite[Theorem 2.9.7]{Federer}, for
$\sigma$-almost every $p$,
\begin{equation}\label{eq:density-series}
 \lim_{r\to 0}\frac{(\Delta w)(B_r(p))}{r^{n-2}}
 =\omega_{n-2} f(p).
\end{equation}
Indeed, $\sigma(B_r(p))/(\omega_{n-2}r^{n-2})\to1$, and the
density of $|F|\,dx$ relative to $\sigma$ is zero because
$\sigma$ is supported on the Lebesgue-null set $\C$ and $|F| \in L^1_{\rm loc}(B_{3/4})$.

We compute the limit in \eqref{eq:density-series} by taking limit of a radial
average. Choose a smooth radial function $\varphi(y)=\varphi_0(|y|)$
supported in $B_1$ with $\int\varphi=1$, and put
\begin{equation}\label{eq.def-ar}
 a(r)=\int w(p+ry)\varphi(y)\,dy = r^{-n} \int w(x)\varphi((x-p)/r)\,dx
\end{equation}
and
\begin{equation*}
    \Psi(y)=\int_{|y|}^1 t\varphi_0(t)\,dt, \qquad |y| \le 1.
\end{equation*}
Extend $\Psi$ by zero for $|y| \ge 1$. It is smooth and
$\nabla\Psi=-y\varphi$, so
$\Delta\Psi=-n\varphi-y\cdot\nabla\varphi$. As a result, 
\begin{equation*}
    \frac{d}{dr} \big( r^{-n} \varphi(x/r) \big) = r^{1-n} \Delta_x \big( \Psi(x/r) \big).
\end{equation*}
Using this, we compute $a'(r)$ from \eqref{eq.def-ar} and apply Fubini's theorem to the integral against the measure $\Delta w$ to get
\begin{equation}
\begin{aligned}
 ra'(r)&=
 r^{2-n} \int w (x) \Delta_x \left\{
 \Psi \left(\frac{x-p}{r}\right) \right\} dx\\
 &=r^{2-n}\int \Psi\left(\frac{x-p}{r}\right)d(\Delta w)(x)
 \\
 &=\int_0^1 t\varphi_0(t)
      \frac{(\Delta w)(B_{rt}(p))}{r^{n-2}}\,dt.
 \label{eq:radial-derivative}
\end{aligned}
\end{equation}
At a point where \eqref{eq:density-series} holds, the quotient in
\eqref{eq:radial-derivative} is bounded by
$C_p t^{n-2}$ for small $r$. Dominated convergence theorem and
$\int\varphi = 1$ yield
\begin{equation}\label{eq.a'limit}
 \lim_{r\to 0}ra'(r)
 =\frac{\omega_{n-2} f(p)}{|\mathbb S^{n-1}|} = (2\pi)^{-1} f(p).
\end{equation}
Here $|\mathbb S^{n-1}|$ arises as the result of $\int_0^1 t^{n-1} \varphi_0(t)\,dt$, due to $\int \varphi_0(|y|) \, dy = 1$.

On the other hand, \eqref{eq:local-growth-log} and \eqref{eq.def-ar} together give $a(r)=k\log r+O_p(1)$. Hence, by L'H\^{o}pital's rule and \eqref{eq.a'limit}
\begin{equation}
    k = \lim_{r\to 0} \frac{a(r)}{\log r} = \lim_{r\to 0} r a'(r) = \frac{f(p)}{2\pi}.
\end{equation}
This proves $f(p) = 2\pi k = 2\pi (m(p) - 1)$ for $\Haus^{n-2}$ almost every $p\in \C$ and hence $\mu \restriction_{\C} = 2\pi (m - 1) \Haus^{n-2}\restriction_\C$. 
Inserting this into \eqref{eq:measure-decomposition}, we obtain \eqref{eq:distribution} on all interior balls, and hence on $\Omega$.
\end{proof}

\begin{remark}
    As a consequence of \eqref{eq:reg}, local $L^1$ convergence of $w_\varepsilon$ and dominated convergence in the term containing $F$ give
    \begin{equation}\label{eq:reglimit}
    \rho_\varepsilon\,dx = \frac{\varepsilon^2|\nabla^2u|_{\HS}^2}
    {(g^2+\varepsilon^2)^2}\, dx \to \mu \restriction_\C = 2\pi(m-    1)\,\Haus^{n-2}\restriction_{\C},
    \end{equation}
    in the sense of distribution.
\end{remark}

\begin{remark}\label{rmk.F}
    The estimate \eqref{eq.F+L1} provides the $L^1$ bound of $(\Delta \log |\nabla u|)_+$. Together with the bound of negative part in  Corollary \ref{cor:q}, we obtain
    \begin{equation}
        \int_{B_{1/2}\setminus \C} |\Delta \log |\nabla u|| = \int_{B_{1/2} \setminus \C} |F|  \le C_n N^2 \log^2(N+2).
    \end{equation}
\end{remark}

\begin{proof}[Proof of Theorem \ref{thm:main}]
    Take $0\leq\chi\leq1$ with $\chi\in C_0^\infty(B_{3/4})$ and $\chi=1$ on $B_{1/2}$. Testing \eqref{eq:distribution} against $\chi$ and using Lemma \ref{lem:sublevel} and Corollary \ref{cor:q} yield
    \begin{equation}
    \begin{aligned}
    2\pi\int_{\C\cap B_{1/2}}(m-1)\,d\Haus^{n-2}
    &\leq \int_{B_{3/4}} \log |g| \Delta\chi\,dx-\int_{B_{3/4}} \chi F\,dx \\
    &\leq C_n N + \int_{B_{3/4} \setminus \C} q\,dx \\
    & \leq C_n N^2 \log^2(N+2).
    \label{eq:men-bound}
    \end{aligned}
    \end{equation}
    This ends the proof.
\end{proof}

\bibliography{ref}

@book {LL01,
    AUTHOR = {Lieb, Elliott H. and Loss, Michael},
     TITLE = {Analysis},
    SERIES = {Graduate Studies in Mathematics},
    VOLUME = {14},
   EDITION = {Second},
 PUBLISHER = {American Mathematical Society, Providence, RI},
      YEAR = {2001},
     PAGES = {xxii+346},
      ISBN = {0-8218-2783-9},
   MRCLASS = {00A05 (26-01 28-01 31-01 35J10 42-01)},
  MRNUMBER = {1817225},
       DOI = {10.1090/gsm/014},
       URL = {https://doi.org/10.1090/gsm/014},
}

@article {Alessandrini,
    AUTHOR = {Alessandrini, Giovanni},
     TITLE = {Critical points of solutions of elliptic equations in two
              variables},
   JOURNAL = {Ann. Scuola Norm. Sup. Pisa Cl. Sci. (4)},
  FJOURNAL = {Annali della Scuola Normale Superiore di Pisa. Classe di
              Scienze. Serie IV},
    VOLUME = {14},
      YEAR = {1987},
    NUMBER = {2},
     PAGES = {229--256 (1988)},
      ISSN = {0391-173X,2036-2145},
   MRCLASS = {35J25 (35B05)},
  MRNUMBER = {939628},
MRREVIEWER = {C.\ A.\ Swanson},
       URL = {http://www.numdam.org/item?id=ASNSP_1987_4_14_2_229_0},
}

@article {AM92,
    AUTHOR = {Alessandrini, G. and Magnanini, R.},
     TITLE = {The index of isolated critical points and solutions of
              elliptic equations in the plane},
   JOURNAL = {Ann. Scuola Norm. Sup. Pisa Cl. Sci. (4)},
  FJOURNAL = {Annali della Scuola Normale Superiore di Pisa. Classe di
              Scienze. Serie IV},
    VOLUME = {19},
      YEAR = {1992},
    NUMBER = {4},
     PAGES = {567--589},
      ISSN = {0391-173X,2036-2145},
   MRCLASS = {35J25 (35B05 35B30 58E05)},
  MRNUMBER = {1205884},
MRREVIEWER = {Xin\ Kang\ Guo},
       URL = {http://www.numdam.org/item?id=ASNSP_1992_4_19_4_567_0},
}

@article {NV17,
    AUTHOR = {Naber, Aaron and Valtorta, Daniele},
     TITLE = {Volume estimates on the critical sets of solutions to elliptic
              {PDE}s},
   JOURNAL = {Comm. Pure Appl. Math.},
  FJOURNAL = {Communications on Pure and Applied Mathematics},
    VOLUME = {70},
      YEAR = {2017},
    NUMBER = {10},
     PAGES = {1835--1897},
      ISSN = {0010-3640,1097-0312},
   MRCLASS = {35J15 (35A20 35R01)},
  MRNUMBER = {3688031},
MRREVIEWER = {Xinhua\ Ji},
       DOI = {10.1002/cpa.21708},
       URL = {https://doi.org/10.1002/cpa.21708},
}

@article {CNV15,
    AUTHOR = {Cheeger, Jeff and Naber, Aaron and Valtorta, Daniele},
     TITLE = {Critical sets of elliptic equations},
   JOURNAL = {Comm. Pure Appl. Math.},
  FJOURNAL = {Communications on Pure and Applied Mathematics},
    VOLUME = {68},
      YEAR = {2015},
    NUMBER = {2},
     PAGES = {173--209},
      ISSN = {0010-3640,1097-0312},
   MRCLASS = {35J25 (35B05)},
  MRNUMBER = {3298662},
MRREVIEWER = {Marco\ Bramanti},
       DOI = {10.1002/cpa.21518},
       URL = {https://doi.org/10.1002/cpa.21518},
}

@article {Dong,
    AUTHOR = {Dong, Rui-Tao},
     TITLE = {Nodal sets of eigenfunctions on {R}iemann surfaces},
   JOURNAL = {J. Differential Geom.},
  FJOURNAL = {Journal of Differential Geometry},
    VOLUME = {36},
      YEAR = {1992},
    NUMBER = {2},
     PAGES = {493--506},
      ISSN = {0022-040X,1945-743X},
   MRCLASS = {58G25 (35P99)},
  MRNUMBER = {1180391},
MRREVIEWER = {Stig\ I.\ Andersson},
       URL = {http://projecteuclid.org/euclid.jdg/1214448750},
}

@book {Federer,
    AUTHOR = {Federer, Herbert},
     TITLE = {Geometric measure theory},
    SERIES = {Die Grundlehren der mathematischen Wissenschaften, Band 153},
 PUBLISHER = {Springer-Verlag New York, Inc., New York},
      YEAR = {1969},
     PAGES = {xiv+676},
   MRCLASS = {28.80 (26.00)},
  MRNUMBER = {257325},
MRREVIEWER = {J.\ E.\ Brothers},
}

@article {Fos25,
    AUTHOR = {Foster, Benjamin},
     TITLE = {Quantitative properties of gradients of solutions to elliptic PDEs},
   JOURNAL = {Ph.D. dissertation, Stanford University},
      YEAR = {2025},
       URL = {},
}

@article {HHHN,
    AUTHOR = {Hardt, R. and Hoffmann-Ostenhof, M. and Hoffmann-Ostenhof, T.
              and Nadirashvili, N.},
     TITLE = {Critical sets of solutions to elliptic equations},
   JOURNAL = {J. Differential Geom.},
  FJOURNAL = {Journal of Differential Geometry},
    VOLUME = {51},
      YEAR = {1999},
    NUMBER = {2},
     PAGES = {359--373},
      ISSN = {0022-040X,1945-743X},
   MRCLASS = {35J15 (35B65)},
  MRNUMBER = {1728303},
MRREVIEWER = {Rolando\ Magnanini},
       URL = {http://projecteuclid.org/euclid.jdg/1214425070},
}

@article {HS89,
    AUTHOR = {Hardt, Robert and Simon, Leon},
     TITLE = {Nodal sets for solutions of elliptic equations},
   JOURNAL = {J. Differential Geom.},
  FJOURNAL = {Journal of Differential Geometry},
    VOLUME = {30},
      YEAR = {1989},
    NUMBER = {2},
     PAGES = {505--522},
      ISSN = {0022-040X,1945-743X},
   MRCLASS = {58E05 (35J99)},
  MRNUMBER = {1010169},
MRREVIEWER = {Fang\ Hua\ Lin},
       URL = {http://projecteuclid.org/euclid.jdg/1214443599},
}

@article {Lin,
    AUTHOR = {Lin, Fang-Hua},
     TITLE = {Nodal sets of solutions of elliptic and parabolic equations},
   JOURNAL = {Comm. Pure Appl. Math.},
  FJOURNAL = {Communications on Pure and Applied Mathematics},
    VOLUME = {44},
      YEAR = {1991},
    NUMBER = {3},
     PAGES = {287--308},
      ISSN = {0010-3640,1097-0312},
   MRCLASS = {58G11 (35J05 35K05 58G03)},
  MRNUMBER = {1090434},
MRREVIEWER = {Robert\ McOwen},
       DOI = {10.1002/cpa.3160440303},
       URL = {https://doi.org/10.1002/cpa.3160440303},
}

@inproceedings {LM18,
    AUTHOR = {Logunov, Alexander and Malinnikova, Eugenia},
     TITLE = {Quantitative propagation of smallness for solutions of
              elliptic equations},
 BOOKTITLE = {Proceedings of the {I}nternational {C}ongress of
              {M}athematicians---{R}io de {J}aneiro 2018. {V}ol. {III}.
              {I}nvited lectures},
     PAGES = {2391--2411},
 PUBLISHER = {World Sci. Publ., Hackensack, NJ},
      YEAR = {2018},
      ISBN = {978-981-3272-92-7; 978-981-3272-87-3},
   MRCLASS = {58J05},
  MRNUMBER = {3966855},
}

@incollection {Sot,
    AUTHOR = {Sottile, Frank},
     TITLE = {From enumerative geometry to solving systems of polynomial equations},
 BOOKTITLE = {Computations in algebraic geometry with {M}acaulay 2},
    SERIES = {Algorithms Comput. Math.},
    VOLUME = {8},
     PAGES = {101--129},
 PUBLISHER = {Springer, Berlin},
      YEAR = {2002},
      ISBN = {3-540-42230-7},
   MRCLASS = {14-04 (13-04 14N10 14N15)},
  MRNUMBER = {1949550},
       DOI = {10.1007/978-3-662-04851-1\{_}6}

@article {Han07,
    AUTHOR = {Han, Qing},
     TITLE = {Nodal sets of harmonic functions},
   JOURNAL = {Pure Appl. Math. Q.},
  FJOURNAL = {Pure and Applied Mathematics Quarterly},
    VOLUME = {3},
      YEAR = {2007},
    NUMBER = {3},
     PAGES = {647--688},
      ISSN = {1558-8599,1558-8602},
   MRCLASS = {31B05},
  MRNUMBER = {2351641},
MRREVIEWER = {Stephen\ J.\ Gardiner},
       DOI = {10.4310/PAMQ.2007.v3.n3.a2},
       URL = {https://doi.org/10.4310/PAMQ.2007.v3.n3.a2},
}

@book{HLbook,
    author    = {Qing Han and Fanghua Lin},
    title     = {Nodal Sets of Solutions of Elliptic Equations},
    publisher = {Unpublished book},
    year      = {2013},
}

@book{SimonGMT,
    author    = {Simon, Leon},
    title     = {Introduction to Geometric Measure Theory},
    publisher = {Lecture notes},
    year      = {2017-2018},
    URL = {https://math.stanford.edu/~lms/ntu-gmt-text.pdf},
}

@article {CM01,
    AUTHOR = {Colding, Tobias H. and Minicozzi, II, William P.},
     TITLE = {Minimal annuli with and without slits},
   JOURNAL = {J. Symplectic Geom.},
  FJOURNAL = {The Journal of Symplectic Geometry},
    VOLUME = {1},
      YEAR = {2001},
    NUMBER = {1},
     PAGES = {47--61},
      ISSN = {1527-5256,1540-2347},
   MRCLASS = {53A10},
  MRNUMBER = {1959578},
MRREVIEWER = {Zhongmin\ Shen},
       DOI = {10.4310/jsg.2001.v1.n1.a2},
       URL = {https://doi.org/10.4310/jsg.2001.v1.n1.a2},
}

@article {CZZ25,
    AUTHOR = {Cai, Yingying and Zhu, Jiuyi and Zhuge, Jinping},
     TITLE = {Quantitative unique continuation for {N}eumann problem in planar ${C}^{1,\alpha}$ domains},
   JOURNAL = {arXiv:2512.07297},
      YEAR = {2025},
       URL = {},
}
\bibliographystyle{abbrv}
\end{document}